\documentclass{amsart}
\usepackage{amssymb,latexsym}
\numberwithin{equation}{section}

\newtheorem{theorem}{Theorem}[section]
\newtheorem{lemma}[theorem]{Lemma}
\newtheorem{proposition}[theorem]{Proposition}
\newtheorem{corollary}[theorem]{Corollary}

\theoremstyle{definition}
\newtheorem{definition}[theorem]{Definition}
\newtheorem{example}[theorem]{Example}

\theoremstyle{remark}
\newtheorem{remark}[theorem]{Remark}

\def\tr{{Trop}}
\def\R{{\mathcal R}}
\def\uq{U_q(\g)}\def\uqn{U_q(\n)}\def\uqnm{U_q(\nm)}
\def\uqo{U_q^0}
\def\uqb{U_q(\b)}\def\uqbm{U_q(\b^-)}
\def\vq{V_q}\def\vql{\vq(\lambda)}

\def\N{{\mathbb Z}_{\geq 0}}
\def\Z{{\mathbb Z}}%
\def\Q{{\mathbb Q}}
\def\n{{\mathfrak n}}%
\def\C{{\mathbb C}}%
\def\g{{\mathfrak g}}%
\def\h{{\mathfrak h}}%
\def\b{{\mathfrak b}}\def\nm{{\mathfrak n}^-}
\def\li{{\mathcal L}_{{\bf i}}}\def\ci{c_{{\bf i}}}\def\Ci{{\mathcal C}_{{\bf i}}}
\def\bi{b_{{\bf i}}}\def\Ais{A_{{\bf i}}^\sigma}
\def\B{{\mathcal B}}\def\E{{\mathcal E}}

\def\K{{\mathbb K}}

\begin{document}
\title[Realisation of Lusztig cones]
{Realisation of Lusztig cones}
\author{Philippe Caldero}
\address{D\'epartement de Math\'ematiques, Universit\'e Claude Bernard Lyon I,
69622 Villeurbanne Cedex, France}
\email{caldero@igd.univ-lyon1.fr}
\author{Robert Marsh}
\address{Department of Mathematics, University of
Leicester, University Road, Leicester LE1 7RH, England}
\email{R.Marsh@mcs.le.ac.uk}
\author{Sophie Morier-Genoud}
\address{D\'epartement de Math\'ematiques, Universit\'e Claude Bernard Lyon I,
69622 Villeurbanne Cedex, France}
\email{morier@igd.univ-lyon1.fr}

\begin{abstract}
Let $U_q(\mathfrak{g})$ be the quantised enveloping algebra associated
to a simple Lie algebra $\mathfrak{g}$ over $\mathbb{C}$. The negative part
$U^-$ of $U_q(\mathfrak{g})$ possesses a canonical basis $\mathcal{B}$ with
favourable properties. Lusztig has associated a cone~\cite{lusztig2}
to a reduced expression $\mathbf{i}$ for the longest
element $w_0$ in the Weyl group of $\mathfrak{g}$, with good properties
with respect to monomial elements of $\mathcal{B}$.
The first author~\cite{caldero2} has associated a subalgebra
$A_{\mathbf{i}}$ of $U^-$, compatible with the dual basis
$\mathcal{B}^*$, to each reduced expression $\mathbf{i}$. We show that, after
a certain twisting, the string parametrisation of the adapted basis
of this subalgebra coincides with the corresponding Lusztig cone. As an
application, we give explicit expressions for the generators of the Lusztig
cones.
\end{abstract}

\date{\today}

\maketitle
\begin{section} {Introduction}
Let $U=U_q({\mathfrak{g}})$ be the quantum group associated to a semisimple
Lie algebra $\mathfrak{g}$. The negative part $U^-$ of $U$ has
a canonical basis $\mathcal{B}$ with favourable properties
(see Kashiwara~\cite{kashiwara2} and Lusztig~\cite[\S14.4.6]{lusztig1}).
For example, via action on highest weight vectors it gives rise to bases
for all the finite-dimensional irreducible highest weight $U$-modules.
The dual canonical basis $\mathcal{B}^*$ of the positive part $U^+$ has good multiplicative properties.
Two elements of $\mathcal{B}^*$ are said to be multiplicative if their
product also lies in $\mathcal{B}^*$ up to a power of $q$. \par
The first author has shown that for each reduced expression
$\mathbf{i}=(i_1,i_2,\ldots ,i_N)$ for the longest element $w_0$
(see \S~\ref{lusztigconedefinition}) in the Weyl group of $\mathfrak{g}$,
there is a corresponding subalgebra $A_{\mathbf{i}}$ of $U^+$, known as a
standard
adapted subalgebra, with basis given by $A_{\mathbf{i}}\cap \mathcal{B}^*$,
consisting entirely of elements which are pairwise multiplicative. The subalgebras $A_{\mathbf{i}}$ are $q$-polynomial algebras, i.e. algebras given by generators and $q$-commuting relations, with GK-dimension $N=l(w_0)$.
Note that adapted algebras were introduced for the Berenstein-Zelevinsky
conjecture and are connected with the larger theory of cluster
algebras~\cite{fominzelevinsky1}.
\par
By a Lusztig cone of $U$, we mean the cone
$\li\subseteq \mathbb{N}^N$ associated by Lusztig
(see~\cite[\S16]{lusztig2}) to each reduced expression
$\mathbf{i}$ for $w_0$.
In~\cite{lusztig2} these cones arise naturally from the linear term of a
nonhomogeneous quadratic form associated to ${\bf i}$ which is used by
Lusztig to give a positivity condition for a monomial
\begin{equation}
F_{i_1}^{(a_1)}F_{i_2}^{(a_2)}\cdots F_{i_N}^{(a_N)}\label{mon}
\end{equation}
to lie in the canonical basis.
Here the $F_i$ are the standard generators of $U^-$. Monomials of this
form with $(a_1,a_2,\ldots ,a_N)$ lying in the Lusztig cone corresponding
to $\mathbf{i}$ lie in the canonical basis
in types $A_1,A_2$ and $A_3$~\cite{lusztig2}, in type
$A_4$~\cite{marsh2} and in type $B_2$~\cite{xi2}; see also~\cite{reineke2}.
Counter-examples of M.~Reineke~\cite{reineke2} and N. H. Xi~\cite{xi3}
show that this fails in general in type $A_5$. Recently,
R.~Bedard~\cite{bedard2} has analysed the quadratic forms associated to these
monomials and, as an application, was able to compute some interesting
examples in types $D_4$, $A_5$ and affine $A_1$.

Note that Lusztig cones are used to describe regular functions on a
reduced real double Bruhat cell of the corresponding algebraic
group~\cite{zelevinsky1}, they have links with primitive elements in the dual
canonical basis (this can be seen using~\cite{berensteinzelevinsky1}) and
therefore with the representation theory of affine Hecke algebras~\cite{lnt1},
and they are known to correspond to regions of linearity of the Lusztig
reparametrisation functions (see~\cite{cartermarsh1}).

Given a reduced expression $\mathbf{i}$ for $w_0$, elements of the dual
canonical basis can be parametrised via the string
parametrisation in direction $\mathbf{i}$
(see~\cite[\S2]{berensteinzelevinsky1},~\cite{kashiwara3}
and~\cite[\S2]{nakashimazelevinsky1}),
which we denote by $c_{\mathbf{i}}:\mathcal{B}^*\rightarrow \mathcal{C}_{\mathbf{i}}$,
where $\mathcal{C}_{\mathbf{i}}\subseteq\mathbb{N}^N$ is known as the string cone
cor\-res\-ponding to $\mathbf{i}$.

Our main result is that the set of string parameters (in direction
$\mathbf{i}$) of a certain twisting of the standard adapted subalgebra of
the dual canonical basis corresponding to $\mathbf{i}$ coincides with the
Lusztig cone corresponding to $\mathbf{i}$. The twisting is done with the help of the Sch\"utzenberger involution.
This gives a realisation of all Lusztig cones in terms of the dual canonical
basis. It also implies that all Lusztig cones are simplicial 
(generalising results of Bedard~\cite{bedard1}
and the second author~\cite{marsh1}), and enables us to give an explicit
description of their spanning vectors; see Theorem~\ref{simplicial}.

The paper is organised as follows. Sections 1 and 2 give preliminary
results on quantum groups and the canonical basis, including its
parametrisations associated to a reduced word, and adapted algebras.
In Section 4, we introduce the Sch\"utzenberger involution $\phi$ and its
action on the dual canonical basis. In Section 5, we recall some facts on
geometric lifting of the canonical basis in order to give a formula which
describes $\phi$ in terms of the parametrisation of the dual canonical basis.
A remarkable property is that, with a good choice of parametrisations, the
action of the Sch\"utzenberger involution on the dual canonical basis is given
by an affine map.

In Sections 6 and 7, we apply the results from previous sections to describe
explicitly the twisted standard adapted subalgebra associated to a reduced
word $\mathbf{i}$, in terms of $\mathbf{i}$-string parametrisation.
By the multiplicative property of the adapted subalgebra and the "affine map"
property, this can be provided by an $N\times N$ matrix and a column vector.
A combinatorial argument, together with the known PBW-parametrisation of
the adapted basis of a standard adapted subalgebra, allows us to prove 
the main theorem: in Section 8, we realise the Lusztig cones in terms of
the string parametrisation of twisted standard adapted subalgebras. As an
application, we give an explicit formula for the generators of the cones. 
\end{section}
\begin{section} {Notation and preliminaries.}
\begin{subsection}{ }
Let $A=(a_{ij})_{1\leq i,j\leq n}$ be the Cartan matrix of a finite
dimensional semi-simple Lie algebra $\g$ over $\C$. Let $\g=\nm\oplus\h\oplus\n$ be a triangular decomposition, where $\h$ is a Cartan subalgebra and where $\nm$, $\n$ are opposite maximal nilpotent subalgebras of $\g$. Let
$\{\alpha_i\}_i$ be the set of simple roots of the root system $\Delta$
resulting from this decomposition.
The set of positive roots is denoted by $\Delta^+$. \par
Let $P$ be the weight lattice generated by the fundamental weights $\varpi_i$, $1\leq i\leq n$. Set $P^+:=\sum_i\N\varpi_i$, endowed with the ordering $\sum_i\lambda_i\varpi_i\leq\sum_i\mu_i\varpi_i\Leftrightarrow
\lambda_i\leq\mu_i$. The Weyl group $W$ is generated by the reflections $s_i$ corresponding to the simple roots. We denote by $<\,,\,>$ the $W$-invariant form on $P$; we have $a_{ij}=<\alpha_j,\alpha_i^{\vee}>$ for all $i,j$.\par
For $n$ a nonnegative integer and $\alpha$ a positive root, we set :
$q_\alpha=q^{<\alpha,\alpha>/2}$,
$[n]_\alpha=\frac{q_\alpha^n-q_\alpha^{-n}}{q_\alpha-q_\alpha^{-1}}$,
$[n]_\alpha!=[n]_\alpha [n-1]_\alpha\ldots[1]_\alpha$.\medskip\noindent
\end{subsection}
\begin{subsection}{ }\label{lusztigconedefinition}

Let $W$ be the Weyl group of $\mathfrak{g}$, with Coxeter generators
$s_1,s_2,\ldots,s_n$ and corresponding length function. An expression
$s_{i_1}s_{i_2}\cdots s_{i_m}$ for an element of $w$ is called reduced if
it is of minimal length; we identify such an expression with the tuple
$\mathbf{i}=(i_1,i_2,\ldots ,i_m)$. Set $N:=\dim\n$. It is known that $N$
is the length of the longest element $w_0$ of the Weyl group. 
Let ${\mathcal R}$ be the set of reduced expressions for $w_0$.

Fix ${\bf i}$ in ${\mathcal R}$. Let $\li$ be the set of points $(c_1,\ldots,c_N)\in\N^N$ with the following property : for any two indices $p<p'$ in $\{1,\ldots,N\}$ such that $i_p=i_{p'}=i$ and $i_q\not= i$ whenever $p<q<p'$, we have
$$c_p+c_{p'}+\sum_{p<q<p'} a_{i_p,i_q}c_q\leq 0.$$
The cone $\li$ is the so-called Lusztig cone associated to the reduced
expression ${\bf i}$. This is defined in~\cite[\S16]{lusztig2} for the
simply-laced case. We use here a natural generalisation to the general case
which also appears implicitly (for type $B_2$) in~\cite{xi2}.
\end{subsection}
\begin{subsection}{ }
Let $d$ be an integer such that $<P,P>\subset (2/d)\Z$.
Let $q$ be a indeterminate and set $\K=\C(q^{1/d})$.
We define the simply connected quantised
enveloping $\K$-algebra $\uq$ as in~\cite{joseph1}. Set
$d_i=<\alpha_i,\alpha_i>/2$
and $q_i=q^{d_i}$ for all $i$. 
Let $\uqn$, resp. $\uqnm$, 
be the subalgebra generated by the canonical generators 
$E_i:=E_{\alpha_i}$, resp. $F_i:=F_{\alpha_i}$,
of positive, resp. negative, weights, subject to the quantum Serre relations.
For all $\lambda$
in $P$, let $K_\lambda$ be the corresponding element in the 
algebra $\uqo=\K[P]$ of the torus of $\uq$ and $K_i:=K_{\alpha_i}$.
We have the triangular decomposition 
$\uq=\uqnm\otimes\uqo\otimes\uqn$. Set 
$$\uqb=\uqn\otimes\uqo, \hskip 15mm \uqbm=\uqnm\otimes\uqo.$$\par
The algebra $\uq$ is endowed with a structure of Hopf algebra with 
comultiplication $\Delta$, 
antipode $S$ and augmentation $\varepsilon$ given by
$$\Delta E_i=E_i\otimes 1+K_i\otimes E_i,\ \Delta F_i
=F_i\otimes K_i^{-1}+ 1\otimes F_i,\ \Delta K_\lambda
=K_\lambda\otimes K_\lambda,$$
$$S(E_i)=-K_i^{-1}E_i,\ S(F_i)=-F_iK_i,\ S(K_\lambda)
=K_{-\lambda},$$
$$\varepsilon(E_i)=\varepsilon(F_i)=0,\,\varepsilon(K_\lambda)
=1.$$\medskip\noindent
Let $(\,,\,)$ be the Hopf bilinear form,~\cite{rosso1}, on $\uqb\times\uqbm$, uniquely defined by 
$$(E_i,F_j)=\delta_{ij}(1-q_i^2)^{-1}, \; 1\leq i,j\leq n,$$
$$(XK_\lambda,YK_\mu)=q^{(\lambda,\mu)}(X,Y), \;X\in\uqn,\,Y\in\uqnm.$$
\end{subsection}
\begin{subsection}{ }\label{automorphisms}
In this section, we define automorphisms of the quantised enveloping algebra
and the Poincar\'e-Birkhoff-Witt basis.
The automomorphisms $T_i$,
$1\leq
i\leq
n$, as in~\cite{lusztig1}, are given by : 
$$T_i(E_i)=-K_{i}^{-1}F_i,\hskip 5mm T_i(E_j)=\sum_{k+l=-a_{ij}}
(-1)^{k}\frac{{q_{\alpha_i}^{-k}}}{[k]_{\alpha_i}![l]_{\alpha_i}!}E_i^{k}E_jE_i^{l},\;
1\leq i,j\leq n,\, i\not=j,$$
$$T_i(F_i)=-E_iK_i,\hskip 3mm T_i(F_j)=
\sum_{k+l=-a_{ij}}
\frac{q_{\alpha_i}^{k}}{[l]_{\alpha_i}![k]_{\alpha_i}!}
F_i^{l}F_jF_i^{k},\; 1\leq i,j\leq n,\,
i\not=j,$$
$$T_i(K_{\alpha_j})=K_{s_i(\alpha_j)},\,1\leq i,j\leq n.$$
It is known that the $T_i$'s define a braid action on $\uq$.
Fix ${\bf i}=(i_1,\ldots,i_N)$ in ${\mathcal R}$. For each $k$, $1\leq k\leq N$, set $\beta_k:=s_{i_1}\ldots
s_{i_{k-1}}(\alpha_k)$. It is well known that $\{\beta_k,\,
1\leq k\leq N\}$ is the set of positive roots and that
$$\beta_1<\beta_2<\ldots<\beta_N$$
defines a so-called convex ordering on $\Delta^+$. This ordering identifies
the semigroup $\N^{\Delta^+}$ with the semigroup $\N^N$. In the sequel, we denote
by $\{e_k, 1\leq k\leq N\}$ the natural basis of this semigroup.\par
For all $k$, define $E^{{\bf i}}_{\beta_k}=E_{\beta_k}=T_{i_1}\ldots
T_{i_{k-1}}(E_{\alpha_{i_k}})$.
For all $t=(t_i)\in\N^N$, set $E^{{\bf i}}(t)=E(t):=E_{\beta_1}^{(t_1)}\ldots
E_{\beta_N}^{(t_N)}$, where $E_{\beta_k}^{(t_k)}:=\frac{1}{[t_k]_{\beta_k}!}E_{\beta_k}^{t_k}$. It is known that $\{E(t), \, t\in\N^N\}$
is a basis of $\uqn$ called the Poincar\'e-Birkhoff-Witt basis, in short
PBW-basis,
associated to the reduced expression ${\bf i}$. In the same way, we
can define the PBW-basis $\{F(t), \, t\in\N^N\}$ of $\uqnm$.
\smallskip
We define now the automorphisms $\bar{\ }\,$, $\omega$, and the antiautomorphism
$\sigma$ of $\uq$ by
$$\bar E_i=E_i, \,\bar K_i=K_i^{-1},\,\bar F_i=F_i, \,\bar q=q^{-1},\, 1\leq i\leq n,$$
$$\omega(E_i)=F_i, \,\omega(K_i)=K_i^{-1},\,\omega(F_i)=E_i, \,\omega(q)=q,\, 1\leq i\leq n,$$
$$\sigma(E_i)=E_i, \,\sigma(K_i)=K_i^{-1},\,\sigma(F_i)=F_i, \,\sigma(q)=q,\, 1\leq i\leq n.$$
Note that $\omega$ is a coalgebra antiautomorphism.
\end{subsection}
\begin{subsection}{}
For any (left) $\uq$-module $M$, and any weight $\mu$ in $P$, let
$M_\mu=\{m\in M\,:\,K_\lambda.m=q^{<\lambda,\mu>}m\}$ be the subspace of $M$
of weight $\mu$. For all $\lambda$ in $P^+$ let $\vql$ be the simple
$\uq$-module with highest weight $\lambda$ and highest weight vector
$v_\lambda$. The module $\vql$ satisfies the Weyl character formula.
For $\lambda$ in $P^+$, $V_q(\lambda)^*$ is naturally endowed with a
structure of right $\uq$-module. Let $\eta_\lambda$ be its weight element
such that $\eta_\lambda(v_\lambda)=1$. For $\lambda$, $\mu$ in $P^+$, the
module $V_q(\lambda)\otimes V_q(\mu)$, endowed with the diagonal action of
$\uq$, has a unique component of type $V_q(\lambda+\mu)$. Hence the
restriction map provides a map $r_{\lambda,\mu}$ :
$V_q(\lambda)^*\otimes V_q(\mu)^*\simeq (V_q(\lambda)\otimes V_q(\mu))^*\rightarrow V_q(\lambda+\mu)^*$.
We can suppose that the $v_\lambda$ are normalised such that
$r_{\lambda,\mu}(\eta_\lambda\otimes\eta_\mu)=\eta_{\lambda+\mu}$,
$\lambda$, $\mu\in P^+$.
\par
Set $R^+:=\bigoplus_{\lambda\in P^+}\vql^*\otimes v_\lambda$. The space $R^+$ can be equipped with a structure of algebra by the multiplication rule : $ (\xi\otimes v_{\lambda}).(\xi'\otimes v_{\lambda'})=r_{\lambda,\lambda'}(\xi\otimes\xi')\otimes v_{\lambda+\lambda'}$, $\xi\in V_q(\lambda)^*$, $\xi'\in V_q(\lambda')^*$. \par
The map $R^+\rightarrow \uqbm^*$, $\xi\otimes v_\lambda\mapsto \xi(?v_\lambda)$, $\xi\in \vql^*$, provides an embedding of algebras.
\par
Define the map 
$$\beta\,:\,\uqb\rightarrow \uqbm^*,\hskip 5mm \beta(u)(v)=(u,v).$$
Then we have:
\begin{theorem}\label{beta} \cite{caldero1}
The map $\beta$ is an injective antihomomorphism of algebras which maps
$K_{\lambda}$ to $\eta_{\lambda}\otimes v_\lambda\in R^+\subseteq \uqbm^*$.
There exists a unique subspace $E_\lambda$ of $\uqn$ such that
$\beta(E_\lambda K_{\lambda})=\vql^*\otimes v_\lambda$.
\end{theorem}
\end{subsection}
\end{section}
\begin{section} {Recollection on canonical bases and adapted algebras.}
The recollection on canonical bases is from~\cite{berensteinzelevinsky2}. 
Results on adapted algebras can be found in~\cite{caldero2}
and~\cite{caldero3}.
\begin{subsection}{ }

The following theorem defines the canonical basis and its Lusztig parametrisation.
\begin{theorem} \cite{lusztig1}
Fix ${\bf i}$ in ${\mathcal R}$. For all $t$ in $\N^N$, there exists a unique element $b=\bi(t)$ in $\uqnm$ such that $\bar b=b$ and $b-F^{{\bf i}}(t)\in
q^{-1}\sum_{t'<t}\Z[q^{-1}]F^{{\bf i}}(t')$. The map $t\mapsto \bi(t)$ defines a bijection from $\N^N$ to a basis $\B$ of $\uqnm$ which does not depend on the choice of ${\bf i}$.
\end{theorem}
The basis $\B$ is called the canonical (or global) basis and the map
$\bi$ : $t\mapsto \bi(t)$ is the Lusztig parametrisation of $\B$ associated to
the reduced expression ${\bf i}$.
We can define the action of the Kashiwara operators on the canonical basis as
follows : for $1\leq i\leq n$, there exists a unique injective map
$\tilde f_i$ :
$\B\rightarrow\B$, such that for all ${\bf i}$ with $i_1=i$, we have
$$\tilde f_i(\bi(t_1,t_2,\ldots,t_N))=\bi(t_1+1,t_2,\ldots,t_N).$$
For $1\leq i\leq n$, let $\tilde e_i$ : $\B\rightarrow\B\cup\{0\}$ be such that
$\tilde e_i(b)=b'$ if there exists $b'$ such that $\tilde f_i(b')=b$ and
$\tilde e_i(b)=0$ if not.\par
Let $\varepsilon_i(b)=$max$\{k,\,\tilde e_i(b)\not=0\}$ and let $\E$ : $\B\rightarrow P^+$, $b\mapsto \sum_i\varepsilon_i(b)\varpi_i$. Now, the basis $\B$ is stable under $\sigma$. For all $\lambda$ in $P^+$, set 
$\B(\lambda)=\{b\in\B,\,\E(\sigma(b))\leq\lambda\}$.
\par
A nice theorem of compatibility of the canonical basis with the Weyl modules $\vql$ can be stated as follows :
\begin{theorem} \cite{kashiwara1}\label{blambda}
Fix $\lambda$ in $P^+$. Then, for $b$ in $\B$, we have $b\cdot v_\lambda\not=0$
if and only if $b\in\B(\lambda)$. Moreover, $\B(\lambda)\cdot v_\lambda$ is a
basis of $\vql$.
\end{theorem}
In the sequel, we will identify $\B(\lambda)$ with its image in $\vql$.
\end{subsection}
\begin{subsection}{ }

We now introduce the string parametrisation of the canonical basis and the various transition maps.\par
Fix a reduced expression ${\bf i}$ in $\R$ and $b$ in $\B$. The string of
$b$ in the direction ${\bf i}$ is the sequence of integers $\ci(b):=(t_1,\ldots,t_N)$ defined recursively by :
$$t_1=\varepsilon_{i_1}(b),\,t_2=\varepsilon_{i_2}(\tilde e_{i_1}^{t_1}(b)),\ldots,
t_N=\varepsilon_{i_N}(\tilde e_{i_{N-1}}^{t_{N-1}}\ldots \tilde e_{i_1}^{t_1}(b)).$$
The map $\ci$ defines a bijection from $\B$ onto the set of integral points of a rational convex polyhedral cone $\Ci$ in $\R^N$.
\par
We can now define :
$$R_{{\bf i}}^{{\bf i}'}=(b_{{\bf i}'})^{-1}\circ\bi\;:\;\N^N\rightarrow\N^N,$$
$$R_{-{\bf i}}^{-{\bf i}'}=c_{{\bf i}'}\circ(\ci)^{-1}\;:\;\Ci\rightarrow{\mathcal C}_{{\bf i}'},$$
$$R_{-{\bf i}}^{{\bf i}'}=(b_{{\bf i}'})^{-1}\circ(\ci)^{-1}\;:\;\Ci\rightarrow\N^N,$$
$$R_{{\bf i}}^{-{\bf i}'}=c_{{\bf i}'}\circ\bi\;:\;\N^N\rightarrow{\mathcal C}_{{\bf i}'}.$$
\end{subsection}
\begin{subsection}{ }

Let $\B^*\subset\uqn$ be the basis dual to $\B$ with respect to the form
$(\,,\,)$ on $\uqn\times\uqnm$. We call it the dual canonical basis.
For $b$ in $\B$, we denote by $b^*$ the corresponding element in $\B^*$.
Since we work with the dual canonical basis, in the sequel we shall regard
$b_{\mathbf{i}}$ as a map from $\mathbb{Z}_{\geq 0}^N$
to $\mathcal{B}^*$ (rather than $\mathcal{B}$), using the identification
$b\leftrightarrow b^*$. Similarly we shall regard $c_{\mathbf{i}}$ as a map
from $\mathcal{B}^*$ to $\mathbb{Z}_{\geq 0}^N$.

The set $\B^*$ is stable under $\sigma$ "up to a power of $q$". To be more precise, for $b$ in $\B^*$, there exists an integer $m$ such that $\sigma(b^*)=q^m\sigma(b)^*.$ \par
For $\lambda$ in $P^+$ and $b$ in $\B(\lambda)$, let $\pi_{\lambda}(b)^*$ be the element of $\vql^*$ such that $\pi_{\lambda}(b)^*(b'.v_\lambda)=\delta_{b,b'}$, where $\delta$ is the Kronecker symbol.
It is easily seen from the definitions that :
\begin{lemma}\label{betacan}
For all $\lambda$ in $P^+$ and $b$ in $\B(\lambda)$, we have $\beta(b^*K_{\lambda})=\pi_{\lambda}(b)^*\otimes v_\lambda$.
\end{lemma}
In the notation, we will sometimes omit $\pi_\lambda$.\par 
The lemma implies that the spaces $E_\lambda$ defined by Theorem~\ref{beta} are compatible with the dual canonical basis. By the Weyl character formula, for all $w$ in $W$ and $\lambda$ in $P^+$, there is a unique element of $\B^*\cap E_\lambda$ with weight $\lambda-w\lambda$. We denote this element by $b_{w,\lambda}^*$ and the corresponding element in the canonical basis by $b_{w,\lambda}$.\par
In the sequel, $\{\pi_\lambda(b)^* \otimes v_\lambda,\,b\in\B(\lambda),\,\lambda\in P^+\}$ will be called dual canonical basis of $R^+$. By a misuse of language, its Lusztig, resp. string, parametrisation will be the Lusztig, resp. string, parametrisation of the corresponding element $b$ in $\B(\lambda)$. 
\end{subsection}
\begin{subsection}{} \label{adaptedalgebras}
Two elements of the dual canonical basis are called multiplicative if their
product is an element of the dual canonical basis, up to a power of $q$.
By~\cite{reineke1}, multiplicative elements $q$-commute.
We start with the definition of adapted algebras.
\begin{definition}
A subalgebra $A$ of $\uqn$, resp. $R^+$, is called {\em adapted}
if \par\noindent
1) the intersection of $A$ and the dual canonical basis is a basis of $A$, called adapted basis, \par\noindent
2) the elements of this basis are pairwise multiplicative.
\end{definition}
We define standard adapted subalgebras of $R^+$, associated to a fixed reduced expression ${\bf i}$ for $w_0$; see~\cite{caldero2}.\par
Set $y_\lambda:=\eta_\lambda\otimes v_\lambda$. Let $A_{{\bf i}}$ be the subalgebra of $R^+$ generated by $y_i:=y_{\varpi_i}$, $1\leq i\leq n$, and $c_k^{{\bf i}}:=b_{s_{i_1}\ldots s_{i_k},\varpi_{i_k}}^*\otimes v_{\varpi_{i_k}}$, $1\leq k\leq N$. (This notation is standard, but shouldn't be confused with
$c_{\mathbf{i}}$ which is used for the string parametrisation).
Using the antihomomorphism $\beta$, we obtain from~\cite{caldero2} the following proposition : 

\begin{proposition}
Fix ${\bf i}\in\R$. The algebra $A_{{\bf i}}$ is an adapted subalgebra of
$R^+$ and the adapted basis is given by monomials in the $y_i$ and the
$c_k^{{\bf i}}$, up to a power of $q$.
\end{proposition}
As $\sigma$ is an antiautomorphism which preserves the dual canonical basis up
to a power of $q$, we can now define another family of adapted algebras by
twisting the standard ones.\par

Note that $\sigma(b_{s_{i_1}\ldots s_{i_k},\varpi_{i_k}})\in\B(\mu_k)$,
where $\mu_k:=\E(b_{s_{i_1}\ldots s_{i_k},\varpi_{i_k}})$.
Let $A_{{\bf i}}^\sigma$ be the subalgebra generated by the $y_i$,
$1\leq i\leq n$, and the
$c_k^{{\bf i}\sigma}:=\sigma(b_{s_{i_1}\ldots s_{i_k},\varpi_{i_k}})^*\otimes v_{\mu_k}$, $1\leq k\leq N$.
These elements $q$-commute by~\ref{beta}.
Moreover, by~\ref{betacan} and~\cite[2.2]{caldero2}, we have
\begin{proposition} \label{sigmaadapted}
Fix ${\bf i}\in\R$. Then,\par
\item{(i)} the algebra $A_{{\bf i}}^\sigma$ is an adapted subalgebra of $R^+$,\par
\item{(ii)} the adapted basis is given by monomials in the $y_i$ and
the $c_k^{{\bf i}\sigma}$, up to a power of $q$,\par
\item{(iii)} the Lusztig parametrisation of $y_i$ is zero,\par
\item{(iv)} the Lusztig parametrisation of $c_k^{{\bf i}\sigma}$ is $\sum e_l$ where $l$ runs over $\{l\leq k,\, i_l=i_k\}$.
\end{proposition}
\end{subsection}
\begin{subsection}{} \label{qcentre}
In this section, we are concerned with the lowest weight vectors in $R^+$.
For all $\lambda$ in $P^+$, set $\lambda^*:=-w_0\lambda$. Let $v_{w_0\lambda}$ be the the unique element of weight $w_0\lambda$ in the canonical basis of $\vql$. Then, $z_\lambda:=v_{w_0\lambda}^*\otimes v_\lambda$ is an element of the dual canonical basis of $R^+$. It is known (see~\cite{caldero2}) that these
elements belong to all of the standard adapted algebras defined in the
section above.
\begin{lemma}\label{qcenter} 
The elements $z_\lambda$, $\lambda\in P^+$ satisfy
$z_\lambda z_\mu=z_{\lambda+\mu}$, $\lambda$, $\mu\in P^+$.
\end{lemma}
\begin{proof}
Let $b$, $b'$ be two elements of the canonical basis of $\uqnm$, with $b\in\B(\lambda)$. By~\cite[Proposition 3.1]{caldero3}, we have $b'^*b^*\in q^{-<\lambda,\nu'>}\B^*$, where $\nu'$ is the weight of $b'^*$.
\par Applying this formula when $b^*$, resp. $b'^*$, is the element of $\B(\lambda)^*$, resp. $\B(\mu)^*$, of weight $\lambda+\lambda^*$, resp. $\mu+\mu^*$, and using Theorem~\ref{beta}, we obtain that
$$z_\mu z_\lambda=q^{-<\lambda,\mu+\mu^*>}q^{<\lambda,\mu+\mu^*>}z_{\lambda+\mu}=z_{\lambda+\mu}.$$
This proves the lemma.
\end{proof}
In the sequel, we set $z_i=z_{\varpi_i}$.
\end{subsection}
\end{section}
\begin{section} {The Sch\"utzenberger involution.}
We can now define an involution $\phi$ of $R^+$ by twisting the dual Weyl modules by the automorphism $\omega$. This involution generalises the Sch\"utzenberger involution,~\cite{schutzenberger1}, up to a diagram automorphism.
\begin{subsection}{}
Let $M$ be a $\uq$-module. With the help of the automorphism $\omega$, we define a twisted $\uq$-module structure on $M^*$ via
$$x.\xi(m)=\xi(\omega(x)m),\,\xi\in M^*,\, x\in\uq,\, m\in M.$$
Let $M_\omega^*$ be the space $M^*$ endowed with this $\uq$-module
structure.\par
Fix $\lambda$ in $P^+$. Then, $\vql_\omega^*$ is a simple right module. By dualising~\cite[Chap XXI]{lusztig1}, we obtain
\begin{proposition}\label{philambda}
The module $\vql_\omega^*$ is isomorphic to $V_q(\lambda^*)^*$. There exists a unique right $\uq$-module isomorphism $\phi_\lambda$ which sends the dual canonical basis of $\vql_\omega^*$ to the dual canonical basis of $V_q(\lambda^*)^*$. It sends highest weight vectors to lowest weight vectors and conversely.
\end{proposition}
\end{subsection}
\begin{subsection}{}
We define a map $\phi$ : $R^+\rightarrow R^+$, such that
$$\phi(\xi\otimes v_\lambda)=\phi_\lambda(\xi)\otimes v_{\lambda^*},\, \xi\in\vql^*,\,\lambda\in P^+,$$
where $\phi_\lambda$ is as above. Then,
\begin{proposition}
The map $\phi$ is an involutive antiautomorphism of the algebra $R^+$ which preserves the dual canonical basis of $R^+$.
\end{proposition}
\begin{proof}
Recall that $\omega$ is involutive. So, by Proposition~\ref{philambda}, $\phi_{\lambda^*}\phi_\lambda$ is the identity and this implies that $\phi$ is involutive.\par
Now, let $R_\lambda^+$ be the $\lambda$-component of $R^+$, which is isomorphic to $V_q(\lambda)^*$ as a right $\uq$-module. As noted in \ref{automorphisms}, $\omega$ is a coalgebra antiautomorphism. Hence, the map $m$ : 
$$ m\,:\,R_\lambda^+\otimes R_\mu^+\rightarrow R_{\lambda+\mu}^+,\;\; m(a\otimes b)=\phi^{-1}(\phi(b)\phi(a))$$ is a morphism of right $\uq$-modules, where $R_\lambda^+\otimes R_\mu^+$ is endowed with the dia\-go\-nal action. Using Lemma~\ref{qcenter}, we obtain $$m(y_\lambda\otimes y_\mu)=\phi^{-1}(z_{\lambda+\mu})=y_{\lambda+\mu}=y_\lambda y_\mu.$$ As $\hbox{dim}\hbox{Hom}_{\uq}(R_\lambda^+\otimes R_\mu^+,R_{\lambda+\mu}^+)=1$, this proves that $m$ is the multiplication of $R^+$ and thus that $\phi$ is an algebra antiautomorphism. 
The last assertion of the proposition is clear by Proposition~\ref{philambda}. 
\end{proof}
\begin{corollary} \label{twistedsubalgebra1}
Fix a reduced expression ${\bf i}$ in $\R$. Then, $\phi(\Ais)$ is an adapted
subalgebra of $R^+$.
\end{corollary}

The aim of the remaining sections is to prove that the Lusztig cone $\li$
is the ${\bf i}$-string parametrisation of the adapted basis of
${\phi}(\Ais)$. We note that this is given by monomials in the $z_i$
and the $\phi(c_k^{\mathbf{i}\sigma})$ up to a power of $q$.
\end{subsection}
\end{section}
\begin{section}{Geometric Lifting and parametrisation.}
Let's fix a dominant weight $\lambda$. In this section, we study the geometric lifting of the morphism $\phi_\lambda$. By using results of~\cite{berensteinzelevinsky2}, we obtain an explicit formula for the morphism $b_{{\bf i}}^{-1}\phi _\lambda c_{{\bf i}}^{-1}$ which gives the Lusztig parametrisation $t'=(t'_1,\cdots, t'_N)$ of the element $\phi _\lambda (b)$ in terms of the string $t=(t_1,\cdots ,t_N)$ of $b$, where $b$ is in the dual canonical basis of $R^+$. 
\subsection{}
We give here notation and recollection of~\cite{berensteinzelevinsky2}.
Let $G$ be the semisimple simply connected complex Lie group with Lie algebra
$\g$. For all $i$, $1\leq i\leq n$, we denote by $\varphi _i:SL_2\hookrightarrow G$ the canonical embedding corresponding to the simple root $\alpha _i$.
Consider the one parameter subgroups of $G$ defined by 
$$x_i(t)=\varphi _i\left(
\begin{array}{ccc}
1&t\\
0&1
\end{array}
\right),
\hspace{0.5cm} y_i(t)=\varphi _i\left(
\begin{array}{ccc}
1&0\\
t&1
\end{array}
\right)
,\hspace{0.5cm}t\in \C,$$

and
$$t^{\alpha _i^\vee}=\varphi _i\left(
\begin{array}{ccc}
t&0\\
0&t^{-1}
\end{array}
\right)
,\hspace{0.5cm}t\in \C^*.$$
The $x_i(t)$, (resp. $y_i(t)$, $t^{\alpha _i ^\vee}$) generate subgroups $N$,
(resp. $N^-$, $H$). 
We have the following commutation relations :
\begin{eqnarray}\label{commutation}
t^{\alpha _i^\vee}x_j(t')=x_j(t^{a_{ij}}t')t^{\alpha _i^\vee}, \hspace{0.3cm}t^{\alpha _i^\vee }y_j(t')=y_j(t^{-a_{ij}}t')t^{\alpha _i^\vee}. 
\end{eqnarray}
We define two involutive antiautomorphisms of $G$, $x\rightarrow x^T$ and $x\rightarrow x^\iota $, by :
$$x_i(t)^T=y_i(t), \hspace{0.3cm}y_i(t)^T=x_i(t), \hspace{0.3cm}(t^{\alpha _i^\vee})^T=t^{\alpha _i^\vee},$$
$$x_i(t)^ \iota =x_i(t), \hspace{0.3cm}y_i(t)^\iota =y_i(t), \hspace{0.3cm}(t^{\alpha _i^\vee})^\iota =t^{-\alpha _i^\vee}.$$
The first one is called transposition and the second one is called inversion.\par
Let $G_0:=N^-HN$ be the set of elements in $G$ which have a (unique) gaussian decomposition; we write $x=[x]_-[x]_0[x]_+$ for the gaussian decomposition of $x$ in $G_0$.

For all reduced expressions ${\bf i}=(i_1,\cdots, i_N) $ and all $N$-tuples $t=(t_1,\cdots, t_N) $ in $\C^N$, we set:
$$x_{{\bf i}}(t):=x_{i_1}(t_1)\cdots x_{i_N}(t_N), \hspace{0.3cm} \text{and} \hspace{0.3cm} 
x_{-{\bf i}}(t):=y_{i_1}(t_1)t_1^{-\alpha _{i_1} ^{\vee}} \cdots y_{i_N}(t_N)
t_N^{-\alpha _{i_N} ^\vee}.$$
The $x_{{\bf i}}$ and the $x_{-{\bf i}}$ parametrise subvarieties of $G$. 
\begin{theorem} \cite{berensteinzelevinsky2}
There exists a subvariety of $G$, $L^{e,w_0}_{>0}$, resp. $L^{w_0,e}_{>0}$, such that for all ${\bf i}$ in $\mathcal{R}$, the map $x_{{\bf i}}$, resp. $x_{-{\bf i}}$, is a bijection from ${\mathbb
R}^N_{>0} $ to $L^{e,w_0}_{>0}$, resp. $L^{w_0,e}_{>0}$.
\end{theorem}

We denote by $\tilde {R}_{{\bf i}}^{{\bf i}'}:=x_{{\bf i}'}^{-1}\circ x_{\bf i}$ and $\tilde {R}_{-{\bf i}}^{-{\bf i}'}:=x_{-{\bf i}'}^{-1}\circ x_{-{\bf i}}$ the transition maps. 
A remarkable result of~\cite{berensteinzelevinsky2} asserts that $\tilde {R}_{{\bf i}}^{{\bf i}'}$ (respectively, $\tilde {R}_{-{\bf i}}^{-{\bf i}'}$) is
a geometric lifting of the map $R_{\bf i}^{{\bf i}'}$ (respectively,
$R_{-{\bf i}}^{-{\bf i}'}$), which was defined in the first section.
Let's be more precise. 
\par
By using results on semifields (see~\cite{berensteinfominzelevinsky1}),
the authors define the so-called tropicalisation, denoted by $[.]_\tr$.
The map $ [.]_\tr$ is a map from the semifield $\Q_{>0}(t_1,\cdots,t_N)$ to
the set of maps $\Z^N\rightarrow \Z$.
The elements of $\Q_{>0}(t_1,\cdots,t_N)$ are called
\textit{subtraction-free rational expressions} in the $t_1,\cdots, t_N$.
Tropicalising a subtraction-free expression means replacing the multiplication
by the operation $a\odot b:=a+b$ and the sum by the operation
$a\oplus b=\min (a,b)$.
Let's give an example from~\cite{berensteinfominzelevinsky1}.
\begin{example}
Let $x$, $y$ be two indeterminates and set $f=x^2-xy+y^2$. Then, $f$ is a
subtraction-free expression because $f=\frac{x^3+y^3}{x+y}$. We have :
$[f]_{\tr}$ : $\Z^2\rightarrow\Z$, with
$$[f]_{\tr}(m,n)=\mbox{Min}(3m,3n) -\mbox{Min}(m,n)=\mbox{Min}(2m,2n).$$
\end{example}
A geometric lifting is an element of the inverse image of this map. We can see in the example above that it is in general not unique. 

The following theorem is a result of~\cite{berensteinzelevinsky2}.
In the theorem, the notation $(.)^\vee$ means that we consider the analogue
formula in the Langlands dual of $G$, and the notation $ [.]_\tr$ is the
componentwise tropicalisation.
\begin{theorem} \label{ThmTropR}
Fix two reduced expressions ${\bf i}$, ${\bf i'}$ in $\R$. Then
$(\tilde {R}_{\bf i}^{{\bf i}'})^\vee$,
resp.
$(\tilde {R}_{-{\bf i}}^{-{\bf i}'})^\vee$,
is a geometric lifting of $R_{\bf i}^{{\bf i}'}$,
resp.
$R_{-{\bf i}}^{-{\bf i}'}$ :

$$ (i)\hspace{0.3cm} [(\tilde {R}_{\bf i}^{{\bf i}'})^\vee (t)]_{\tr}=R_{\bf i}^{{\bf i}'}(t), \hspace{0.5cm} (ii) \hspace{0.3cm} [(\tilde {R}_{-{\bf i}}^{-{\bf i}'})^\vee (t)]_{\tr}=R_{-{\bf i}}^{-{\bf i}'}(t).$$
\end{theorem}

\subsection{}
Let $\zeta :L^{w_0,e}_{>0}\rightarrow L^{e,w_0}_{>0} $ be the map
defined by $$\zeta (x):=[x^{\iota T}]_{+}.$$

By~\ref{commutation}, we obtain that the map $\zeta$ is well defined and 
\begin{proposition} \label{FormZeta} 
Let ${\bf i}=(i_1,\cdots, i_N) $ be a reduced expression for $w_0$, and suppose $(t_1',\cdots ,t_N')\-=(x_{\bf i}^{-1}\circ \zeta \circ x_{-{\bf i}})(t_1,\cdots ,t_N)$, $t_i\in\C$. Then, we have
$$t'_k=t_k^{-1}\prod _{j>k}t_j^{-a_{i_ji_k}}.$$
\end{proposition}

The following theorem and its corollary are a result of~\cite{moriergenoud1}
but we give here a sketch of the proof.
The description of the geometric lifting of $\phi _\lambda $ can be given in
terms of parametrisations :
\begin{theorem}\label{bellaformula}
Fix two reduced expressions ${\bf i}$, ${\bf i'}$ in $\R$.
Then, $(x_{\bf i}^{-1}\circ \zeta \circ x_{-{\bf i}'})^\vee (t)$ is a
subtraction-free expression and
$$ b_{{\bf i}}^{-1}\phi _\lambda c_{{\bf i}'}^{-1}(t)= [(x_{\bf i}^{-1}\circ \zeta \circ x_{-{\bf i}'})^\vee (t)]_\tr+b_{{\bf i}}^{-1}\phi_{\lambda} (\eta_{\lambda}).$$
\end{theorem}

Set $(l_1,\cdots ,l_N):=b_{{\bf i}}^{-1}\phi _\lambda (\eta_\lambda )$. We obtain the following tropicalised formula :
\begin{corollary} \label{bellaformulamatrix}
For $(t_1',\cdots ,t_N')=b_{{\bf i}}^{-1}\phi _\lambda c_{{\bf i}}^{-1}(t_1, \cdots,t_N)$,
$$t'_k=l_k-t_k-\sum _{j>k}a_{i_ki_j}t_j.$$
\end{corollary}
\begin{remark} It is remarkable that this formula is affine. This is only true in the case ${\bf i}={\bf i'}$. In general, the tropical term in the right hand side of \ref{bellaformula} is piecewise linear. 
\end{remark}

\noindent \textit{Sketch of the proof.}
By Proposition~\ref{FormZeta},
$(x_{\bf i}^{-1}\circ \zeta \circ x_{-{\bf i}})$ is a subtraction-free
expression. The first assertion of the theorem is obtained by composing with
$\tilde {R}_{{\bf i}'}^{{\bf i}}$.\par
Let $\phi _{{\bf i},{\bf i}'}:$ ${\mathcal C}_{\bf i'}\rightarrow \Z^N$ be
a family of maps labelled by two reduced expressions ${\bf i}$ and
${\bf i'}$ such that the three following conditions are satisfied :
\begin{enumerate}
\item $\phi _{{\bf i},{\bf i}'}(0,\cdots,0)=b_{{\bf i}}^{-1}\phi _\lambda (\eta_\lambda ),$
\item $\phi _{{\bf i},{\bf i}'}=R_{{\bf i}''}^{\bf i}\circ \phi _{{\bf i}'',{\bf i}'}=\phi _{{\bf i},{\bf i}''}\circ R_{-{\bf i}'}^{-{\bf i}''},$
\item for $\phi _{{\bf i},{\bf i}}(t_1,\cdots,t_N)=(t_1',\cdots ,t_N')$, $t'_1+t_1$ and $t'_k$, $k\not=1$, depend only on $t_2, \cdots, t_N$.
\end{enumerate}
The theorem follows from the proposition:

\begin{proposition} \cite{moriergenoud1}
We have,
\begin{enumerate}
\item[(i)] $(\phi _{{\bf i},{\bf i}'})$ is a family satisfying (1), (2), (3) if and only if
$$ \phi _{{\bf i},{\bf i}'}=b_{{\bf i}}^{-1}\phi _\lambda c_{{\bf i}'}^{-1}.$$
\item [(ii)]The family $(\phi _{{\bf i},{\bf i}'})$ defined by
$$\phi _{{\bf i},{\bf i}'}(t)= [(x_{\bf i}^{-1}\circ \zeta \circ x_{-{\bf i}'})^\vee (t)]_\tr+b_{{\bf i}}^{-1}\phi _\lambda (\eta_\lambda )$$
satisfies the conditions (1), (2), (3).
\end{enumerate}
\end{proposition}
\end{section}
\begin{section}{From the PBW-parametrisation to the string parametrisation}
Let ${\bf i}=(i_1,\cdots, i_N)$ be a reduced expression for $w_0$.
Our aim in this section is to describe the map
$u\mapsto c_{\mathbf{i}}\phi_{\lambda}b_{\mathbf{i}}(u)$, using
Theorem~\ref{bellaformula}.
\par
Let $\mathcal{S}$ be the complex rational map extending
$$(t_1,t_2,\ldots ,t_N)\mapsto (u_1,u_2,\ldots ,u_N)=
(x_{\mathbf{i}}^{-1}\circ \zeta \circ x_{-\mathbf{i}})^{\vee}
(t_1,t_2,\ldots ,t_N)$$
(see Proposition~\ref{FormZeta}). Then Theorem~\ref{bellaformula} states that
$\mathcal{S}(t)$ is a subtraction-free expression
(given by Proposition~\ref{FormZeta}), and that
$$b_{{\bf i}}^{-1}\phi _\lambda c_{{\bf i}}^{-1}(t)=
[\mathcal{S}(t)]_\tr+
b_{{\bf i}}^{-1}\phi_\lambda (\eta_\lambda ).$$

Moreover, $\mathcal{S}$ is clearly birational.
We first explain how to invert $\mathcal{S}$.

\begin{lemma} \label{Sinverse}
Let $(u_1,u_2,\ldots u_N)\in{\C}^N$. Then
$\mathcal{S}^{-1}(u_1,u_2,\ldots u_N)=(t_1,t_2,\ldots t_N)\in\mathbb{C}^N$,
where, for $1\leq k \leq N$, we have
\begin{equation} \label{inverseformula}
t_k=u_k^{-1}\prod_{j>k}u_j^{<s_{i_{j-1}}\cdots s_{i_{k+1}}
\alpha_{i_k},\alpha_{i_j}^{\vee}>}.
\end{equation}
\end{lemma}

\begin{proof}
We note that, by Proposition~\ref{FormZeta}, if
$(t_1,t_2,\ldots ,t_N)\in\mathbb{C}^N$ and $\mathcal{S}(t_1,t_2,\ldots ,t_N)=
(u_1,u_2,\ldots ,u_N)\in\mathbb{C}^N$, then, for $1\leq k\leq N$, we have:
\begin{equation} \label{forwardformula}
u_k=t_k^{-1}\prod _{j>k}t_j^{-a_{i_ji_k}}.
\end{equation}
Since this map is a monomial transformation of ${\C}^N$, it is sufficient to
show that substituting the expression~\eqref{forwardformula} for $u_k$ in
terms of the $t_j$ into right hand side of equation~\eqref{inverseformula}
reduces to the left hand side. We obtain:
$$\left( t_k^{-1}\prod_{l>k} t_l^{-a_{i_l,i_k}} \right)^{-1}
\prod_{j>k}\left( t_j^{-1}\prod_{l>j}t_l^{-a_{i_l,i_j}}\right)^{<s_{i_{j-1}}\cdots s_{i_{k+1}}\alpha_{i_k},\alpha_{i_j}^{\vee}>}.$$
It is clear that the exponent of $t_k$ in this expression is $1$, and that,
if $l<k$, then the exponent of $t_l$ is zero. So we consider the case where
$l>k$. The exponent of $t_l$ is given by
$$
a_{i_l,i_k}+
\left(\sum_{l>j>k}-a_{i_l,i_j}<s_{i_{j-1}}\cdots s_{i_{k+1}}\alpha_{i_k},\alpha_{i_j}^{\vee}>\right)
-<s_{i_{l-1}}\cdots s_{i_{k+1}}\alpha_{i_k},\alpha_{i_l}^{\vee}>.
$$
Since $s_{i_j}(\alpha_{i_l}^{\vee})=\alpha_{i_l}^{\vee}-a_{i_l,i_j}\alpha_{i_j}^{\vee}$, this is equal to:
\begin{eqnarray*}
<\alpha_{i_k},\alpha_{i_l}^{\vee}> &+&
\sum_{l>j>k}<s_{i_{j-1}}\cdots s_{i_{k+1}}\alpha_{i_k},s_{i_j}(\alpha_{i_l}^{\vee})>-<s_{i_{j-1}}\cdots s_{i_{k+1}}\alpha_{i_k},\alpha_{i_l}^{\vee}> \\
&-& <s_{i_{l-1}}\cdots s_{i_{k+1}}\alpha_{i_k},\alpha_{i_l}^{\vee}>.
\end{eqnarray*}
The sum telescopes to give:
$$
<\alpha_{i_k},\alpha_{i_l}^{\vee}>+
<s_{i_{l-1}}\cdots s_{i_{k+1}}\alpha_{i_k},\alpha_{i_l}^{\vee}>
-<\alpha_{i_k},\alpha_{i_l}^{\vee}>
-<s_{i_{l-1}}\cdots s_{i_{k+1}}\alpha_{i_k},\alpha_{i_l}^{\vee}>=0,$$
and we are done.
\end{proof}

\begin{proposition} \label{inversebellaformula}
Fix a reduced expression ${\bf i}$ for $w_0$, and let
$u=(u_1,u_2,\ldots ,u_N)\in \mathbb{Z}_{\geq 0}^N$.
Then
$$c_{\bf i}\phi _\lambda b_{\bf i}(u)=
[(\mathcal{S}^{-1})(u)]_{\tr}+c_{\bf i}\phi_\lambda (\eta_\lambda),$$
where $\mathcal{S}^{-1}(u)$ is as in Lemma~\ref{Sinverse}.
\end{proposition}

\begin{proof}
By Theorem~\ref{bellaformula}, we have that
\begin{equation} \label{quotebella}
b_{{\bf i}}^{-1}\phi _\lambda c_{{\bf i}}^{-1}(t)=
[S(t)]_{trop}+b_{{\bf i}}^{-1}\phi_\lambda (\eta_\lambda).
\end{equation}
It follows that
$$c_{\bf i}\phi _{\lambda^*} b_{\bf i}(u)=[\mathcal{S}^{-1}(u-
b_{\mathbf{i}}^{-1}\phi_{\lambda}(\eta_{\lambda}))]_{\tr},$$
noting that $\phi_{\lambda^*}\phi_{\lambda}$ is the identity map.
We note that $\mathcal{S}$ is an invertible monomial map, so its
tropicalisation is linear.
Hence
$$c_{\bf i}\phi_{\lambda^*} b_{\bf i}(u)=[\mathcal{S}^{-1}(u)]_{\tr}-
[\mathcal{S}^{-1}(b_{\mathbf{i}}^{-1}\phi_{\lambda}(\eta_{\lambda}))]_{\tr}.$$
Substituting $t=c_{\mathbf{i}}\phi_{\lambda^*}(\eta_{\lambda^*})$
into~\eqref{quotebella}, we obtain
$$0=b_{\mathbf{i}}^{-1}\eta_{\lambda^*}=[S(c_{\mathbf{i}}\phi_{\lambda^*}(\eta_{\lambda^*}))]_{\tr}+b_{\mathbf{i}}^{-1}\phi_{\lambda}(\eta_{\lambda}).$$
Hence
$$[\mathcal{S}^{-1}(b_{\mathbf{i}}^{-1}\phi_{\lambda}(\eta_{\lambda}))]_{\tr}=
-c_{\mathbf{i}}\phi_{\lambda^*}(\eta_{\lambda^*}).$$
Hence we have:
$$c_{\bf i}\phi_{\lambda^*} b_{\bf i}(u)=[\mathcal{S}^{-1}(u)]_{\tr}
+c_{\mathbf{i}}\phi_{\lambda^*}(\eta_{\lambda^*}),$$
giving the required result (since $\lambda\mapsto \lambda^*$ is an
involution).
\end{proof}

We can compute the constant term in this formula as follows:
\begin{lemma} \label{loweststring}
Let $\mathbf{i}$ be a reduced expression for $w_0$. Then
$c_{\mathbf{i}}(\phi_{\lambda}\eta_{\lambda})=(v_1,v_2,\ldots ,v_N)$, where,
for $1\leq k\leq N$, we have
$$v_k=<s_{i_{k-1}}\cdots s_{i_1}\lambda,\alpha_{i_k}^{\vee}>.$$
\end{lemma}
\begin{proof}
This follows from~\cite[28.1.4]{lusztig1}, since we are
computing the string of the lowest weight vector in $V_q(\lambda^*)^*$.
\end{proof}
We remark that $\phi_{\lambda}(\eta_{\lambda})=z_{\lambda^*}$, so this
lemma is computing the string $c_{\mathbf{i}}z_{\lambda^*}$.
\end{section}
\begin{section}{The string parametrisation of a twisted standard adapted subalgebra}
Recall (see Corollary~\ref{twistedsubalgebra1})
that the elements $z_1,z_2,\ldots ,z_n$,
together with the elements $\phi(c_k^{\mathbf{i}\sigma})$,
$k=1,2, \ldots ,N$,
form an adapted basis for the twisted standard adapted subalgebra
$\phi A_{\mathbf{i}}^{\sigma}$ of $R^+$. Our aim is to 
compute the string parameters of these elements. We therefore use
Proposition~\ref{inversebellaformula} to apply the
map $c_{\bf i}\phi_{\mu_k} b_{\bf i}$ to the element
$b_{\mathbf{i}}^{-1}(c_k^{\mathbf{i}\sigma})$, for each $1\leq k\leq N$;
see~\ref{adaptedalgebras}.
These last vectors were described in Proposition~\ref{sigmaadapted}(iv).

\begin{subsection}{}
For convenience, we define a matrix $V$ with columns given by these
vectors.

\begin{definition}
Let $M_N(\mathbb{Z})$ denote the ring of $N\times N$ matrices with integer
entries.
Let $V=(V_{jk})\in M_N(\mathbb{Z})$ be defined by:
$$V_{jk}=\left\{
\begin{array}{ll}
1 & j\leq k, \mbox{\ if\ }i_j=i_k, \\
0 & \mbox{otherwise}.
\end{array}
\right.
$$
\end{definition}

Let $\mathbf{v}_k$ denote the $k$th column of $V$; this is $b_{\mathbf{i}}^{-1}(c_k^{\mathbf{i}\sigma})$ by Proposition~\ref{sigmaadapted}(iv).
\end{subsection}
\begin{subsection}{}
In order to apply Proposition~\ref{inversebellaformula} to compute
$c_{\bf i}\phi_{\mu_k} b_{\bf i}(\mathbf{v}_k)$, for each $1\leq k\leq N$, we
first need to apply $(\mathcal{S}^{-1})_{\tr}$ to $\mathbf{v}_k$.

Let $S$ be the matrix defining the linear map $\mathcal{S}_{\tr}$.
Then the matrix of $(\mathcal{S}^{-1})_{\tr}$ is $S^{-1}$.
By Lemma~\ref{Sinverse}, we have that $S^{-1}=T=(T_{jk})$,
where
$$
T_{jk}=
\left\{
\begin{array}{cc}
-1 & j=k, \\
<s_{i_{j+1}}\cdots s_{i_{k-1}}\alpha_{i_k},\alpha_{i_j}^{\vee}> & j<k, \\
0 & \mbox{otherwise.}
\end{array}
\right.
$$

\begin{definition} \label{matrixC}
Let $C=(C_{jk})\in M_N(\mathbb{Z})$ be the matrix given by:
$$
C_{jk}=
\left\{
\begin{array}{cc}
<s_{i_{j+1}}\cdots s_{i_k}\varpi_{i_k},\alpha_{i_j}^{\vee}> & j\leq k, \\
0 & \mbox{otherwise.}
\end{array}.
\right.
$$
\end{definition}

We have:
\begin{lemma} \label{SVC}
Let $S,V$ and $C$ be the matrices defined above. Then we have
$S^{-1}V=-C$.
\end{lemma}
\begin{proof}
We have that $(S^{-1}V)_{jl}=\sum_{j\leq k\leq l} T_{jk}V_{kl}$, and is
zero if $j>l$. If $j\leq l$ then, using that $<\,,\,>$ is $W$-invariant,
\begin{equation}\label{smv}(S^{-1}V)_{jl}=T_{jj}V_{jl}+\sum_{j<k\leq l,i_k=i_l}<s_{i_{k-1}}\cdots
s_{i_{j+1}}\alpha_{i_j}^{\vee},\alpha_{i_l}>.
\end{equation}
Remark that $T_{jj}V_{jl}=-\delta_{i_j,i_l}$. Now, by calculating explicitly the coefficient of $\alpha_{i_l}^{\vee}$ in $s_{i_l}\cdots s_{i_j}\alpha_{i_j}^{\vee}$, we find that it is equal to the right hand side of \ref{smv}. Hence, 
$$(S^{-1}V)_{jl}=<\varpi_{i_l},s_{i_l}\cdots s_{i_j}\alpha_{i_j}^{\vee}>=
-<s_{i_{j+1}}\cdots s_{i_l}\varpi_{i_l},\alpha_{i_j}^{\vee}>,$$
which is equal to $-C_{jl}$ as required.
\end{proof}

We have the following Corollary:

\begin{corollary} \label{Cpositive}
The entries of $C$ are nonnegative.
\end{corollary}
\begin{proof}
In proof of Lemma~\ref{SVC}, we noted that $C_{jl}$ is equal to the negative
of the coefficient of $\alpha_{i_l}^{\vee}$ in the negative coroot
$s_{i_l}\cdots s_{i_j}\alpha_{i_j}^{\vee}$.
\end{proof}
\end{subsection}
\begin{subsection}{}
We now would like to compute $c_{\bf i}\phi _{\mu_k} b_{\bf i}(\mathbf{v}_k)$
for each $k$ (note that the $\mathbf{v}_k$ are the columns of $V$).
In the following lemma we set
$\varepsilon_l(c_k^{\mathbf{i}})=\varepsilon_l(b_{s_{i_1}\ldots s_{i_k},\varpi_{i_k}})$ by a misuse of notation.
\begin{lemma} \label{epsilonvalue}
For $1\leq l\leq n$ and $1\leq k\leq N$, we have:
$$\varepsilon_l(c_k^{\mathbf{i}})=\left\{
\begin{array}{cc}
-<s_{i_1}\cdots s_{i_k}\varpi_{i_k},\alpha_l^{\vee}> &
\mbox{if\ }<s_{i_1}\cdots s_{i_k}\varpi_{i_k},\alpha_l^{\vee}>\leq 0, \\
0 & \mbox{otherwise.}
\end{array}\right.$$
In particular,
$$\mu_k=\sum_{1\leq l\leq n,
<s_{i_1}\cdots s_{i_k}\varpi_{i_k},\alpha_l^{\vee}>\leq 0}
-<s_{i_1}\cdots s_{i_k}\varpi_{i_k},\alpha_l^{\vee}>\varpi_l.$$

\end{lemma}
\begin{proof}
This follows from the definition of the $c_k^{\mathbf{i}}$
and~\cite[28.1.4]{lusztig1} (and $sl_2$-representation theory).
\end{proof}

Let $\mathbf{c}_k$ be the $k$th column of $C=-S^{-1}V$, and
let $\mathbf{p}_k=c_{\mathbf{i}}(\phi_{\mu_k}\eta_{\mu_k})$.
Let $P$ be the matrix with columns $\mathbf{p}_k$, $k=1,2,\ldots ,N$.

\begin{lemma} \label{translation}
For $1\leq k\leq N$, we have
$$c_{\mathbf{i}}\phi_{\mu_k}b_{\mathbf{i}}(\mathbf{v}_k)=
-\mathbf{c}_k+\mathbf{p_k}.$$
In particular, the entries of $-C+P$ are nonnegative.
\end{lemma}
\begin{proof}
This follows immediately from Proposition~\ref{inversebellaformula} and Lemma \ref{SVC}.
\end{proof}

As a consequence, we have:

\begin{proposition} \label{parameters1}
Let $\mathbf{i}$ be any reduced expression for $w_0$.
Then $c_{\mathbf{i}}\phi(A_{\mathbf{i}}^{\sigma})$
coincides with the nonnegative integer span of the columns of $-C+P$
and the strings $c_{\mathbf{i}}(z_i)$, $1\leq i\leq n$.
\end{proposition}
\begin{proof}
By Proposition~\ref{sigmaadapted},
$c_{\mathbf{i}}\phi(A_{\mathbf{i}}^{\sigma})$ is the nonnegative integer span
of the $c_{\mathbf{i}}\phi_{\mu_k}b_{\mathbf{i}}(\mathbf{v}_k)$ together
with the strings $c_{\mathbf{i}}(z_i)$, $1\leq i\leq n$. 
Since $c_{\mathbf{i}} (b^*b'^*)=c_{\mathbf{i}}(b^*)+c_{\mathbf{i}}(b'^*)$
when the elements $b^*$ and $b'^*$ of the dual canonical basis are
multiplicative~\cite[Cor. 3.3]{berensteinzelevinsky1},
the proposition follows from
Lemma~\ref{translation} and Corollary~\ref{twistedsubalgebra1}.
\end{proof}
\end{subsection}
\begin{subsection}{}
We also note the following formula for the entries of $P$:

\begin{lemma} \label{pentries}
For $1\leq j,k\leq N$, we have
$$P_{jk}=\sum_{1\leq l\leq n, 
<s_{i_1}\cdots s_{i_k}\varpi_{i_k},\alpha_l^{\vee}>\leq 0}
-<s_{i_1}\cdots s_{i_k}\varpi_{i_k},\alpha_l^{\vee}>
<s_{i_{j-1}}\cdots s_{i_1}\varpi_l,\alpha_{i_j}^{\vee}>.$$
\end{lemma}

\begin{proof}
By definition $P_{jk}$ is the $j$th entry of $c_{\mathbf i}(\phi_{\mu_k}\eta_{\mu_k})$. Hence, the lemma results from Lemma~\ref{epsilonvalue} and Lemma~\ref{loweststring}.
\end{proof}
\end{subsection}
\begin{subsection}{}
We will next show that some of the columns of $-C+P$ are entirely zero,
and therefore can be neglected in Proposition~\ref{parameters1}.
\begin{definition}
Given $k\in\{1,2,\ldots ,N\}$, we set $k(1)=\min\{j\,:\,j>k,i_j=i_k\}$, i.e.
the first occurrence of $i_k$ to the right of $i_k$ in $\mathbf{i}$.
If there is no such occurrence, we set $k(1)=N+1$.
\end{definition}

\begin{lemma} \label{lastoccurrence}
Suppose that $k(1)=N+1$.
Then, for $j=1,2,\ldots ,N$, we have $P_{jk}=C_{jk}$, i.e.
the $k$th column of $P$ coincides with the $k$th column of $C$.
\end{lemma}
\begin{proof}
We have
$$C_{jk}=
\left\{
\begin{array}{cc}
<s_{i_{j+1}}\cdots s_{i_k}\varpi_{i_k},\alpha_{i_j}^{\vee}> & j\leq k, \\
0 & j>k.
\end{array}
\right.
$$
Since $k(1)=N+1$, we have
$$C_{jk}=<s_{i_{j+1}}\cdots s_{i_N}\varpi_{i_k},\alpha_{i_j}^{\vee}>,$$
(in either case).
By 
Lemma~\ref{epsilonvalue},
\begin{eqnarray*}
\mu_k
& = &
\sum_{1\leq l\leq n,
<s_{i_1}\cdots s_{i_k}\varpi_{i_k},\alpha_l^{\vee}>\leq 0}
-<s_{i_1}\cdots s_{i_k}\varpi_{i_k}\alpha_l^{\vee}>\varpi_l \\
& = &
\sum_{1\leq l\leq n,
<s_{i_1}\cdots s_{i_N}\varpi_{i_k},\alpha_l^{\vee}>\leq 0}
-<s_{i_1}\cdots s_{i_N}\varpi_{i_k}\alpha_l^{\vee}>\varpi_l \\
& = &
\sum_{1\leq l\leq n,
<w_0\varpi_{i_k},\alpha_l^{\vee}>\leq 0}
-<w_0\varpi_{i_k},\alpha_l^{\vee}>\varpi_l \\
& = & -w_0\varpi_{i_k}.
\end{eqnarray*}
Hence, by Lemma \ref{loweststring}
\begin{eqnarray*}
P_{jk} & = & -<s_{i_{j-1}}\cdots s_{i_1}w_0\varpi_{i_k},\alpha_{i_j}^{\vee}> \\
& = & -<s_{i_j}\cdots s_{i_N}\varpi_{i_k},\alpha_{i_j}^{\vee}> \\
& = & <s_{i_{j+1}}\cdots s_{i_N}\varpi_{i_k},\alpha_{i_j}^{\vee}> \\
& = & C_{jk}.
\end{eqnarray*}
\end{proof}
\end{subsection}
\begin{subsection}{}
We replace the zero columns in $-C+P$ with vectors which we will see are
the strings of the elements $z_1,z_2,\ldots ,z_n$, in order to obtain a
matrix whose nonnegative integer span is the set of string parameters of the
twisted standard adapted subalgebra of $R^+$ corresponding to $\mathbf{i}$.
We call this matrix $X$:

\begin{definition} \label{spanningvectors}
Let $X=(X_{jk})\in M_N(\mathbb{Z})$ be the matrix defined as follows:
$$X_{jk}=\left\{
\begin{array}{cc}
<s_{i_{j-1}}\cdots s_{i_1}\varpi_{i_k},\alpha_{i_j}^{\vee}>, & k(1)=N+1, \\
-C_{jk}+P_{jk}, & \mbox{otherwise.}
\end{array}
\right.
$$
We note that the $P_{jk}$ are given by Lemma~\ref{pentries}, and that
the $C_{jk}$ are given in Definition~\ref{matrixC}. Also, it follows from
Lemma~\ref{loweststring} that if $k(1)=N+1$ then the $k$th column of
$X$ is the string $c_{\mathbf{i}}z_{\varpi_{i_k}^*}$.
\end{definition}

We have:

\begin{proposition} \label{parameters2}
Let $\mathbf{i}$ be any reduced expression for $w_0$.
Then $c_{\mathbf{i}}{\phi}(A_{\mathbf{i}}^{\sigma})$
coincides with the nonnegative integer span of the columns of $X$.
In particular, the entries of $X$ are nonnegative.
\end{proposition}
\begin{proof}
We note that the matrix $X$ is the same as $-C+P$, except that if
$1\leq k\leq N$ and $k(1)=N+1$, then the $k$th column (which is zero by
Lemma~\ref{lastoccurrence}) is replaced by the string of
$z_{\varpi_{i_k}^*}$.
The result now follows from Proposition~\ref{parameters1}.
\end{proof}
\end{subsection}
\end{section}
\begin{section}{Lusztig cones and twisted standard adapted subalgebras}

In this section, we will show that the cone of string parameters of
$c_{\mathbf{i}}{\phi}(A_{\mathbf{i}}^{\sigma})$, given by the
nonnegative integer span of the columns of $X$, coincides with the Lusztig cone
corresponding to $\mathbf{i}$. At the same time we will show that the
Lusztig cones are simplicial.

\begin{subsection}{}
We first of all define a matrix $\widetilde{L}$
whose rows include the defining inequalities of the Lusztig cone
corresponding to $\mathbf{i}$ as a subset of $\mathbb{N}^N$. This matrix
will later be modified to a matrix defining the Lusztig cone as a subset
of $\mathbb{Z}^N$.

\begin{definition}
Let ${\widetilde{L}}\in M_N(\mathbb{Z})$ be the matrix
defined by:
$$
\widetilde{L}_{jk}=
\left\{
\begin{array}{cc}
-1 & k=j \mbox{\ or\ } k=j(1) \\
-a_{i_j,i_k} & j<k<j(1), \\
0 & \mbox{otherwise.}
\end{array}
\right.
$$
\end{definition}

\begin{remark}
Let $\widetilde{\mathbf{r}}_j$ denote the $j$th row of $\widetilde{L}$.
The defining inequalities of the Lusztig cone
$\mathcal{L}_{\mathbf{i}}$
(as a subset of $\mathbb{Z}_{\geq 0}^N$)
are those inequalities of the form
$\widetilde{\mathbf{r}}_j\cdot \mathbf{c}\geq 0$ for those $j$ such that
$j(1)\leq N$.
\end{remark}
\end{subsection}
\begin{subsection}{}
We will next show how this matrix is related to the matrices $S$ and $V$
already considered, i.e. that $V^{-1}S=\widetilde{L}$. This will have the
consequence that $\widetilde{L}S^{-1}V=I$, in particular showing that the
columns of $S^{-1}V$ satisfy the defining inequalities of the Lusztig cone
corresponding to $\mathbf{i}$.
We recall from Proposition~\ref{FormZeta} that
$\mathcal{S}_{\tr}$ is defined by the matrix $S=(S_{jk})\in M_N(\mathbb{Z})$
where
$$
S_{jk}=
\left\{
\begin{array}{cc}
-1 & j=k \\
-a_{i_j,i_k} & j<k, \\
0 & \mbox{otherwise.}
\end{array}
\right.
$$
We next need to compute the inverse of the matrix $V$.

\begin{lemma} \label{vinverse}
Let $W=(W_{jk})\in M_N(\mathbb{Z})$ be the matrix defined as follows:
$$W_{jk}=\left\{
\begin{array}{cc}
1 & j=k \\
-1 & j<k, k=j(1), \\
0 & \mbox{otherwise.}
\end{array}
\right.
$$
Then $W=V^{-1}$.
\end{lemma}
\begin{proof}
We show that $WV=I$, the identity matrix.
The $j,l$-entry of $WV$ is given by $Z_{jl}=\sum_{k=1}^N W_{jk}V_{kl}$.
For this to be non-zero, we must have $j\leq k\leq l$ and $k=j$ or
$k=j(1)$. There are $5$ cases: \\
Case (a): If $l<j$, then clearly $Z_{jl}=0$. \\
Case (b): If $l=j$, then $Z_{jl}=Z_{jj}=W_{jj}V_{jj}=1\cdot 1=1$. \\
Case (c): If $j<l<j(1)$, then $Z_{jl}=W_{jj}V_{jl}=1\cdot 0=0$. \\
Case (d): If $l=j(1)$, then $Z_{jl}=W_{jj}V_{j,j(1)}+W_{j,j(1)}V_{j(1),j(1)}=
1\cdot 1+(-1)\cdot 1=0$. \\
Case (e): If $l>j(1)$, then $Z_{jl}=W_{jj}V_{jl}+W_{j,j(1)}V_{j(1),l}=
V_{jl}-V_{j(1),l}=0$ since $V_{jl}=V_{j(1),l}$.
\end{proof}

\begin{lemma} \label{VSL}
Let $V,S$ and $\widetilde{L}$ be the matrices as defined above. Then
$V^{-1}S=\widetilde{L}$.
\end{lemma}
\begin{proof}
The $j,l$-entry of $V^{-1}S$ is given by $Y_{jl}=\sum_{k=1}^N W_{jk}S_{kl}$.
To be non-zero, we must have $j\leq k\leq l$ and $k=j$ or $k=j(1)$.
As before, we have the $5$ cases: \\
Case (a): If $l<j$ then clearly $Y_{jl}=0$. \\
Case (b): If $l=j$, then $Y_{jl}=Y_{jj}=W_{jj}S_{jj}=1\cdot (-1)=-1$. \\
Case (c): If $j<l<j(1)$, then $Y_{jl}=W_{jj}S_{jl}=S_{jl}=-a_{i_j,i_l}$. \\
Case (d): If $l=j(1)$, then $Y_{jl}=Y_{j,j(1)}=W_{jj}S_{j,j(1)}+W_{j,j(1)}
S_{j(1),j(1)}=1\cdot (-2)+(-1)\cdot (-1)=-1$. \\
Case (e): If $l>j(1)$, then $Y_{jl}=W_{jj}S_{jl}+W_{j,j(1)}S_{j(1),l}=
S_{jl}-S_{j(1),l}=-a_{i_j,i_l}+a_{i_{j(1)},i_l}=-a_{i_j,i_l}+a_{i_j,i_l}=0$.
\\
We see that $Y_{jl}=\widetilde{L}_{jl}$ in every case.
\end{proof}
\end{subsection}
\begin{subsection}{}
We now define a slightly altered version of the matrix $\widetilde{L}$,
whose rows will eventually be seen to be the defining inequalities of the
Lusztig cone as a subset of $\mathbb{Z}^N$. We will also see that this matrix
is the inverse of the matrix $X$.

\begin{definition}
Let $L=(L_{jk})\in M_N(\mathbb{Z})$ be the matrix defined as follows:
$$
L_{jk}=\left\{
\begin{array}{cc}
-1 & k=j \mbox{\ or\ } k=j(1), \\
-a_{i_j,i_k} & j<k<j(1), \\
1 & j(1)=N+1,\ s_{i_1}s_{i_2}\cdots s_{i_{k-1}}(\alpha_{i_k})=\alpha_{i_j}, \\
0 & \mbox{otherwise.}
\end{array}
\right.
$$
\end{definition}

Firstly, we show how some rows of $L$ are related to the strings of
lowest weight vectors:

\begin{lemma} \label{stringequalities}
Let $$\mathbf{v}=c_{\mathbf{i}}(\phi_{\lambda}\eta_{\lambda})=
(v_1,v_2,\ldots ,v_N),$$ as in Lemma \ref{loweststring}.
Suppose that $1\leq j\leq j(1)\leq N$.
Let $\mathbf{r}_j$ be the $j$th row of $L$.
Then we have $\mathbf{r}_j\cdot \mathbf{v}=0$.
\end{lemma}
\begin{proof}
Recall that for $k=1,2,\ldots ,N$, we have
$$v_k=<s_{i_{k-1}}\cdots s_{i_1}\lambda,\alpha_{i_k}^{\vee}>,$$
by Lemma~\ref{loweststring}.
We have
\begin{eqnarray*}
\mathbf{r}_j\cdot \mathbf{v} & = &
-v_j-v_{j_{(1)}}-\sum_{j<k<j_{(1)}} a_{i_j,i_k}v_k \\
& = & -<s_{i_{j-1}}\cdots s_{i_1}\lambda,\alpha_{i_j}^{\vee}>-
<s_{i_{j(1)-1}}\cdots s_{i_1}\lambda,\alpha_{i_j}^{\vee}>
-\sum_{j<k<j(1)} a_{i_j,i_k}<s_{i_{k-1}}\cdots s_{i_1}\lambda,\alpha_{i_k}^{\vee}>.
\end{eqnarray*}
We note that $$s_{i_k}(\alpha_{i_j}^{\vee})=\alpha_{i_j}^{\vee}-
<\alpha_{i_k},\alpha_{i_j}^{\vee}>\alpha_{i_k}^{\vee}=
\alpha_{i_j}^{\vee}-a_{i_j,i_k}\alpha_{i_k}^{\vee},$$
so
\begin{eqnarray*}
\mathbf{r}_j\cdot \mathbf{v}
& =
& -<s_{i_{j-1}}\cdots s_{i_1}\lambda,\alpha_{i_j}^{\vee}>
-<s_{i_{j(1)-1}}\cdots s_{i_1}\lambda,\alpha_{i_j}^{\vee}> \\
&& +
\sum_{j<k<j(1)}<s_{i_{k-1}}\cdots s_{i_1}\lambda,s_{i_k}(\alpha_{i_j}^{\vee})-
\alpha_{i_j}^{\vee}> \\
& =
& -<s_{i_{j-1}}\cdots s_{i_1}\lambda,\alpha_{i_j}^{\vee}>
-<s_{i_{j(1)-1}}\cdots s_{i_1}\lambda,\alpha_{i_j}^{\vee}> \\
&& +
<s_{i_{j(1)-1}}\cdots s_{i_1}\lambda,\alpha_{i_j}^{\vee}>
-<s_{i_j}\cdots s_{i_1}\lambda,\alpha_{i_j}^{\vee}>=0,
\end{eqnarray*}
the sum telescoping.
\end{proof}

Next, we show that some entries of $-C+P=S^{-1}V+P$ are zero:

\begin{lemma} \label{simplepc}
Suppose that $1\leq j\leq N$ and that
$s_{i_1}s_{i_2}\cdots s_{i_{j-1}}(\alpha_{i_j}^{\vee})=\alpha_r^{\vee}$
is a simple coroot. Then
$$P_{jk}=C_{jk}.$$
\end{lemma}

\begin{proof}
We note that, by Lemma~\ref{pentries},
$$P_{jk}=\left\{
\begin{array}{cc}
-<s_{i_1}\cdots s_{i_k}\varpi_{i_k},\alpha_r^{\vee}> &
<s_{i_1}\cdots s_{i_k}\varpi_{i_k},\alpha_r^{\vee}>\leq 0, \\
0 & <s_{i_1}\cdots s_{i_k}\varpi_{i_k},\alpha_r^{\vee}>>0,
\end{array}
\right.
$$
using the fact that $<\varpi_l,\alpha_r^{\vee}>=\delta_{lr}$.
Suppose first that $j\leq k$. Then:
\begin{eqnarray*}
<s_{i_1}\cdots s_{i_k}\varpi_{i_k},\alpha_r^{\vee}> & = &
<s_{i_{j+1}}\cdots s_{i_k}\varpi_{i_k},s_{i_j}\cdots s_{i_1}\alpha_r^{\vee}> \\
& = & <s_{i_{j+1}}\cdots s_{i_k}\varpi_{i_k},s_{i_j}\alpha_{i_j}^{\vee}> \\
& = & -<s_{i_{j+1}}\cdots s_{i_k}\varpi_{i_k},\alpha_{i_j}^{\vee}> \\
& = & -C_{jk}\leq 0,
\end{eqnarray*}
by Corollary~\ref{Cpositive}. So $P_{jk}=C_{jk}$.
If $j>k$ then
\begin{eqnarray*}
<s_{i_1}\cdots s_{i_k}\varpi_{i_k},\alpha_r^{\vee}> & = &
<\varpi_{i_k},s_{i_k}\cdots s_{i_1}\alpha_r^{\vee}> \\
& = & <s_{i_{j-1}}\cdots s_{i_{k+1}}\varpi_{i_k},s_{i_{j-1}}\cdots s_{i_{k+1}}s_{i_k}
\cdots s_{i_1}\alpha_r^{\vee}> \\
& = & <s_{i_{j-1}}\cdots s_{i_{k+1}}\varpi_{i_k},\alpha_{i_j}^{\vee}> \\
& = &<\varpi_{i_k},s_{i_{k+1}}\cdots s_{i_{j-1}}\alpha_{i_j}^{\vee}>\,\geq 0,
\end{eqnarray*}
since $s_{i_{k+1}}\cdots s_{i_{j+1}}\alpha_{i_j}^{\vee}$ is a positive
coroot. It follows that $P_{jk}=0=C_{jk}$.
\end{proof}

\begin{remark}
We note that the condition
$$s_{i_1}\cdots s_{i_k}(\alpha_{i_j}^{\vee})=\alpha_r^{\vee}$$
is equivalent to the condition
$$s_{i_1}\cdots s_{i_k}(\alpha_{i_j})=\alpha_r.$$
Secondly, this result shows that
$(S^{-1}V+P)_{jk}=0$ under this assumption, by Lemma~\ref{SVC}.
\end{remark}
\end{subsection}
\begin{subsection}{}
We can now prove the following, as we have all the pieces we need:

\begin{proposition} \label{inverseLusztig}
Let $L$ and $X$ be the matrices defined as above, so that the nonnegative
integer span of the columns of $X$ is the cone of the string parameters of
the twisted standard adapted subalgebra corresponding to $\mathbf{i}$.
Then $LX=I$, the identity matrix.
\end{proposition}

\begin{proof}
Denote by $\mathbf{x}_k$ the $k$th column of $X$.
Suppose first that $1\leq j\leq j(1)\leq N$, and that
$1\leq k\leq k(1)\leq N$.
Then $\mathbf{r}_j$ (the $j$th row of $L$) is the same as the $j$th row
of $\widetilde{L}$.
By Lemma~\ref{VSL}, $\widetilde{L}S^{-1}V=I$. Hence
$\mathbf{r}_j\cdot (-\mathbf{c}_k)=\delta_{jk}$, since $-\mathbf{c}_k$ is
the $k$th column of $S^{-1}V$.
By Lemma~\ref{stringequalities}, $\mathbf{r}_j\cdot \mathbf{p}_k=0$ 
(where $\mathbf{p}_k$ is the $k$th column of $P$).
It follows (from the definition of $X$) that
$\mathbf{r}_j\cdot \mathbf{x}_k=\delta_{jk}$.

If $1\leq j\leq j(1)\leq N$ and $k(1)=N+1$, then
$\mathbf{x}_k=c_{\mathbf{i}}z_{\varpi_{i_k}}^*$. By
Lemma~\ref{stringequalities}, we have that
$\mathbf{r}_j\cdot \mathbf{x}_k=0$ (as required, noting that we must have
$j\not=k$).

If $j(1)=N+1$, then let $1\leq l\leq N$ be defined by
$s_{i_1}\cdots s_{i_{l-1}}\alpha_{i_l}=\alpha_{i_j}$. Then
$L_{jl}=1$ is the only non-zero entry in $\mathbf{r}_j$.
It follows that $\mathbf{r}_j.\mathbf{x}_k=X_{lk}$. \\
Case (a) Suppose that $k(1)\leq N$. Then $X_{lk}=-C_{lk}+P_{lk}=0$ by
Lemma~\ref{simplepc}. \\
Case (b) Suppose that $k(1)=N+1$. Then
\begin{eqnarray*}
X_{lk} & = & <s_{i_{l-1}}\cdots s_{i_1}\varpi_{i_k},\alpha_{i_l}^{\vee}> \\
& = & <\varpi_{i_k},s_{i_1}\cdots s_{i_{l-1}}\alpha_{i_l}^{\vee}> \\
& = & <\varpi_{i_k},\alpha_{i_j}^{\vee}> \\
& = & \delta_{j,k},
\end{eqnarray*}
as required (noting that $j(1)=k(1)=N+1$). \\
The proposition is proved.
\end{proof}

If $\mathbf{c}=(c_1,c_2,\ldots c_N)\in\mathbb{Z}^N$, we write
$\mathbf{c}\geq 0$ to denote $c_k\geq 0$ for $k=1,2,\ldots ,N$.
We have the following consequences.

\begin{theorem} \label{simplicial}
The Lusztig cone $\mathcal{L}_{\mathbf{i}}$ is simplicial, defined by the matrix $L$:
$$\mathcal{L}_{\mathbf{i}}=\{\mathbf{c}\in\mathbb{Z}^N\,:\,L\mathbf{c}
\geq 0\}.$$
It coincides with the nonnegative integer span of the columns of the
matrix X (see Definition~\ref{spanningvectors}).
\end{theorem}

\begin{proof}
By Proposition~\ref{parameters2}, the entries of $X$ are nonnegative.
By Proposition~\ref{inverseLusztig}, $XL=I$, so nonnegative integer
combinations of the rows of $L$ are of the form
$(0,0,\ldots ,0,1,0,\ldots ,0)$ (with a $1$ in the $k$th position) and
therefore correspond to inequalities of the form $c_k\geq 0$.
So
$$\{\mathbf{c}\in\mathbb{Z}^N\,:\,L\mathbf{c}\geq 0\}\subseteq \mathbb{N}^N.$$
Since the inequalities corresponding to rows of $L$ are either defining
inequalities of $\mathcal{L}_{\mathbf{i}}$ or inequalities of the form
$c_k\geq 0$, the claimed equality follows, and it is then immediate that
$\mathcal{L}_{\mathbf{i}}$ is spanned by the columns of $L^{-1}=X$.
\end{proof}

\begin{remark}
The fact that $\mathcal{L}_{\mathbf{i}}$ is simplicial was already known for
quiver-compatible reduced expressions for $w_0$ for $\mathbf{g}$ simply
laced~\cite{bedard1} and for all reduced expressions for $w_0$ in type
$A_n$~\cite{marsh1}.
\end{remark}

And we have:

\begin{theorem} \label{parameters3}
Let $\mathbf{i}$ be any reduced expression for $w_0$. Let
$\mathcal{L}_{\mathbf{i}}$ denote the Lusztig cone corresponding to
$\mathbf{i}$.
Let ${\phi}(A_{\mathbf{i}}^{\sigma})$ denote the twisted standard
adapted subalgebra corresponding to $\mathbf{i}$. Then
$$c_{\mathbf{i}}({\phi}(A_{\mathbf{i}}^{\sigma}))=
\mathcal{L}_{\mathbf{i}}.$$
\end{theorem}
\begin{proof}
This follows from Proposition~\ref{parameters2} and Theorem~\ref{simplicial}.
\end{proof}
\end{subsection}
\end{section}

\begin{example}
Suppose that $\mathfrak{g}=sl_4(\mathbb{C})$ (type $A_3$). Let
$\mathbf{i}=(2,3,2,1,2,3)$, a reduced expression for $w_0$. Then we have:
$$\begin{array}{ll}
V=\left(\begin{array}{cccccc}
1 & 0 & 1 & 0 & 1 & 0 \\
0 & 1 & 0 & 0 & 0 & 1 \\
0 & 0 & 1 & 0 & 1 & 0 \\
0 & 0 & 0 & 1 & 0 & 0 \\
0 & 0 & 0 & 0 & 1 & 0 \\
0 & 0 & 0 & 0 & 0 & 1
\end{array}\right),
&
T=
\left(\begin{array}{cccccc}
-1 & -1 & 1 & 0 & -1 & 1 \\
0 & -1 & -1 & -1 & 0 & 1 \\
0 & 0 & -1 & -1 & 1 & 0 \\
0 & 0 & 0 & -1 & -1 & -1 \\
0 & 0 & 0 & 0 & -1 & -1 \\
0 & 0 & 0 & 0 & 0 & -1
\end{array}\right),
\end{array}
$$
$$\begin{array}{ll}
C=-TV=\left(\begin{array}{cccccc}
1 & 1 & 0 & 0 & 1 & 0 \\
0 & 1 & 1 & 1 & 1 & 0 \\
0 & 0 & 1 & 1 & 0 & 0 \\
0 & 0 & 0 & 1 & 1 & 1 \\
0 & 0 & 0 & 0 & 1 & 1 \\
0 & 0 & 0 & 0 & 0 & 1
\end{array}\right),
&
P=
\left(\begin{array}{cccccc}
1 & 1 & 0 & 0 & 1 & 0 \\
1 & 1 & 1 & 1 & 1 & 0 \\
0 & 0 & 1 & 1 & 0 & 0 \\
1 & 1 & 1 & 1 & 1 & 1 \\
1 & 1 & 0 & 0 & 1 & 1 \\
0 & 0 & 0 & 0 & 0 & 1
\end{array}\right).
\end{array}
$$
The matrix $X$ of spanning vectors of the Lusztig cone corresponding
to $\mathbf{i}$ and the defining matrix $L$ are given by:
$$\begin{array}{ll}
X=\left(\begin{array}{cccccc}
0 & 0 & 0 & 0 & 1 & 0 \\
1 & 0 & 0 & 0 & 1 & 1 \\
0 & 0 & 0 & 0 & 0 & 1 \\
1 & 1 & 1 & 1 & 1 & 1 \\
1 & 1 & 0 & 1 & 1 & 0 \\
0 & 0 & 0 & 1 & 0 & 0
\end{array}\right),
&
L=X^{-1}=
\left(\begin{array}{cccccc}
-1 & 1 & -1 & 0 & 0 & 0 \\
0 & -1 & 1 & 0 & 1 & -1 \\
0 & 0 & -1 & 1 & -1 & 0 \\
0 & 0 & 0 & 0 & 0 & 1 \\
1 & 0 & 0 & 0 & 0 & 0 \\
0 & 0 & 1 & 0 & 0 & 0
\end{array}\right).
\end{array}
$$
\end{example}

\end{document}